\documentclass[10pt,a4paper]{extarticle}
\usepackage[utf8]{inputenc}
\usepackage[T1]{fontenc}
\usepackage{ulem}
\usepackage[left=3.5cm,right=3cm,top=2.5cm,bottom=2.5cm]{geometry}
\usepackage{lmodern,mparhack,relsize}     
\usepackage[english]{babel}
\usepackage{amsfonts,amsmath,amsthm,amscd,amssymb,latexsym,mathtools}
\usepackage{epsfig}
\usepackage{ragged2e}
\usepackage{varioref}
\usepackage{caption}
\usepackage{mfirstuc}
\usepackage{subcaption}
\usepackage{xcolor}   
\usepackage[hyphens]{url}
\usepackage[affil-it]{authblk}
\usepackage{hyperref}
\hypersetup{
	breaklinks=true, %
	colorlinks=true, 
	frenchlinks=true,%
	linkcolor=red, %
	citecolor=red, %
	filecolor=magenta, %
	urlcolor=magenta, %
	linktoc=all,     
}
\usepackage{orcidlink} 			
\usepackage{graphicx}
\graphicspath{ {./images/} }
\usepackage{enumerate}
\usepackage{float}
\usepackage[final]{microtype} 
\usepackage[%
	backend=biber,%
    citestyle=numeric-comp,%
	style=numeric-comp,%
	sorting=nyt,%
	hyperref,%
	block=space,%
	backref%
]{biblatex}
\usepackage[autostyle=true]{csquotes}
\usepackage{etoolbox}

\newtheorem{theorem}{Theorem}
\numberwithin{theorem}{section} 

\newtheorem{remark}{Remark}
\numberwithin{remark}{section} 

\newtheorem{definition}{Definition}
\numberwithin{definition}{section} 

\numberwithin{equation}{section}
\numberwithin{figure}{section}

\newcommand{\keyphrase}[1]{\textcolor{blue}{#1}}
\newcommand{\rrspace}[1]{\mathcal{L}\left(#1\right)}
\newcommand{\gfspace}[1]{\mathcal{L}\left(#1\right)}

\newcommand{\ord}[1]{\mathcal{O}\left(#1\right)}
\newcommand{\cost}[1]{\operatorname{C}\left(#1\right)}
\newcommand{\K}{\mathbb{K}}
\newcommand{\GF}[1]{\mathrm{GF}\roundbrack{#1}}
\newcommand{\edwards}{\mathcal{E}}
\newcommand{\edwardsof}[1]{\edwards\roundbrack{#1}}
\newcommand{\weier}{\mathcal{W}}
\newcommand{\weierof}[1]{\weier\roundbrack{#1}}

\newcommand{\gen}[1]{\mathrm{#1}}
\newcommand{\nvector}[1]{\mathbf{#1}}
\newcommand{\roundbrack}[1]{\left(#1\right)}
\newcommand{\divp}[1]{\operatorname{div}\left(#1\right)}
\newcommand{\secref}[1]{\hyperref[sec:#1]{section~\ref{sec:#1}}}
\newcommand{\subsecref}[1]{\hyperref[subsec:#1]{subsec.~\ref{subsec:#1}}}
\newcommand{\figref}[1]{\hyperref[fig:#1]{figure~\ref{fig:#1}}}
\newcommand{\tabref}[1]{\hyperref[tab:#1]{table~\ref{tab:#1}}}
\newcommand{\defref}[1]{\hyperref[def:#1]{definition~\ref{def:#1}}}
\newcommand{\thmref}[1]{\hyperref[thm:#1]{theorem~\ref{thm:#1}}}
\newcommand{\lemmaref}[1]{\hyperref[lemma:#1]{lemma~\ref{lemma:#1}}}
\newcommand{\corollaryref}[1]{\hyperref[corollary:#1]{corollary~\ref{corollary:#1}}}
\newcommand{\remarkref}[1]{\hyperref[rmk:#1]{remark~\ref{rmk:#1}}}
\newcommand*{\email}[1]{\normalsize\href{mailto:#1}{#1}\par}

\everymath{\displaystyle}

\makeatletter
\def\@maketitle{%
	\newpage
	\null
	\vskip 2em%
	\begin{center}%
		\let \footnote \thanks
		{\Large\bfseries \@title \par}%
		\vskip 1.5em%
		{\normalsize
			\lineskip .5em%
			\begin{tabular}[t]{c}%
				\@author
			\end{tabular}\par}%
		\vskip 1em%
		{\normalsize \@date}%
	\end{center}%
	\par}
\makeatother

\title{\bfseries \Large Goppa codes over Edwards curves}
\author{\large Giuseppe Filippone \orcidlink{0000-0001-7315-1852} \thanks{Electronic address: \email{giuseppe.filippone01@unipa.it}}}
\affil{%
	\normalsize Department of Mathematics and Computer Science,\\%
	\normalsize University of Palermo, Via Archirafi 34, 90123 Palermo, Italy}
\date{}

\begin{document}

\maketitle

\begin{abstract}
	Given an Edwards curve, we determine a basis for the Riemann-Roch space of any divisor
	whose support does not contain any of the two singular points.
	This basis allows us to compute a generating matrix for an algebraic-geometric Goppa code over the Edwards curve.
\end{abstract}


\section{Introduction}\label{introduction}
The literature on elliptic curves and their applications in \keyphrase{cryptography} is well consolidated.
Besides the well-known ECC (Elliptic Curve Cryptography) in which the group law, defined on
these curves, is exploited to encrypt messages, and the ECDSA
(Elliptic Curve Digital Signature Algorithm), another example can be found
in the Lenstra algorithm for the factorization of integers.
Moreover, there are as well applications to
\keyphrase{coding theory} based on the \keyphrase{Riemann-Roch space} $ \rrspace{D} $
associated with a rational divisor $ D $ of these curves.
In particular, this space is a fundamental ingredient to construct \keyphrase{Goppa codes},
first introduced in $ 1983 $ \cite{G2}.
Goppa codes over the Hermitian curve, as well as over maximal curves and hyperelliptic curves,
have been extensively studied in \cite{KorchmarosNagyTimpanella, KorchmarosSpeziali, CastellanosFanali, FanaliGiulietti, falcone2020explicit, BartoliBoniniGiulietti, Giulietti2008}, as they have become an important topic both
in coding theory and in cryptography, where they play a central role in \keyphrase{McEliece public-key cryptographic systems}
\cite{Mc}.

To the best of our knowledge, AG Goppa codes for Edwards curves have not been considered until now.
In this paper, we compute the generating matrices for AG Goppa codes
over Edwards curves. These curves are already the subject of many papers in
cryptography \cite{Lange2011, Edwards, BernsteinLange, BernsteinLange2, HisilWongCarterDawson, BernsteinBirknerLangePeters}, in particular in their \keyphrase{twisted}
version. Compared to the classic elliptic curves in Weierstrass form, they can be more efficient
for cryptographic use and for the (single or multiple) digital signature.

In \secref{1} we describe Edwards curves and their relationship with elliptic
curves in Weierstrass form. In \secref{2} we compute a basis
for $ \gfspace{D} $ over Edwards curves, while in \secref{3} we construct AG Goppa codes over Edwards
curves and their generating matrices. In particular,
in \subsecref{goppa_edwards} we give a small example of a
\keyphrase{Goppa MDS code} where we use the AG Goppa code defined in
\subsecref{notations} over Edwards curves.

\section{Edwards curves and elliptic curves in Weierstrass form}\label{sec:1}

In this section we introduce Edwards curves $ \edwards $, that is, algebraic curves,
defined over a field $ \K $, which can be represented in a suitable coordinate system
by the equation $ {\hat{x}}^2 + {\hat{y}}^2 = 1 + d {\hat{x}}^2 {\hat{y}}^2 $,
with $ d (d - 1) \ne 0 $.
We present these curves as a birationally equivalent version of elliptic curves in Weierstrass form.

Recall that, over a field $ \K $ of characteristic different from $ 2 $, a (smooth) elliptic curve
(possessing at least a $ \K $-rational point) can be represented in a suitable coordinate
system by the Weierstrass equation $ y^2 = x^3 + a^\prime x^2 + b^\prime x $, having one point
at infinity $ \Omega = [Z : X : Y] = [0: 0: 1] $ on the $ y $ axis.

\begin{remark}
	Note that, unlike those in Weierstrass form, curves in Edwards form
	$ \edwards $ have two points at infinity, that is,
	$ \Omega_1 = \mathmbox{[\hat{Z} : \hat{X} : \hat{Y}]} = \mathmbox{[0: 1: 0]} $ on
	the $ x $ axis and $ \Omega_2 = \mathmbox{[\hat{Z} : \hat{X} : \hat{Y}]} = \mathmbox{[0: 0: 1]} $
	on the $ y $ axis, which are \keyphrase{ordinary singular points} for
	$ \edwards $	as this curve is \keyphrase{non-smooth}.
\end{remark}

\begin{remark}\label{rmk:remarkable}
	Edwards curves have four remarkable points: $ O = (0, 1) $, $ O^\prime = (0, -1) $,
	$ H = (1, 0) $, $ H^\prime = (-1, 0) $.
	In particular, $ O - O $ is the identity element of the group law defined on them, that is, $ (P - O) + (O - O) = P - O $,
	for any point $ P = (a, b) \in \edwardsof{\K} $, and these four points form the cyclic group
	$ C_4 $, where $ O^\prime - O $ have order $ 2 $, while $ 2 (H - O) = 2 (H^\prime - O) = O^\prime - O $.
\end{remark}

Edwards curves and elliptic curves in Weierstrass form are closely related.
In particular, over a field $ \K $ of characteristic different from $ 2 $,
one has that an elliptic curve $ \weier $ defined by the equation
$ y^2 = x^3 + a^\prime x^2 + b^\prime x $, and an Edwards curve $ \edwards $
defined by the equation $ {\hat{x}}^2 + {\hat{y}}^2 = 1 + d {\hat{x}}^2 {\hat{y}}^2 $,
where $ d $ is not a square, are \keyphrase{birationally equivalent} (cf.\cite{BernsteinLange}).
Furthermore, this equivalence is given by the following two rational maps:
\begin{subequations}
	\begin{align}
		\alpha\colon \edwardsof{\K} &\longrightarrow \weierof{\K} \label{eq:alpha} \\ \nonumber
		(\hat{x}, \hat{y}) &\longmapsto (x, y) = \left(x_1 \dfrac{1 + \hat{y}}{1 - \hat{y}}, y_1 \dfrac{(1 + \hat{y})}{\hat{x} (1 - \hat{y})}\right)\\
		\beta \colon \weierof{\K} &\longrightarrow \edwardsof{\K} \label{eq:beta} \\ \nonumber
		(x, y) &\longmapsto (\hat{x}, \hat{y}) = \left(\dfrac{y_1 x}{x_1 y}, \dfrac{x - x_1}{x + x_1}\right),
	\end{align}
\end{subequations}
where $ P = (x_1, y_1) \in \weierof{\K} $ is such that the divisor
$ 2(P - \Omega) = (0, 0) - \Omega $.

\begin{remark}\label{rmk:birational-equiv-E-W}
	The two rational maps $ \alpha $ and $ \beta $ defines
	a \keyphrase{birational equivalence} between $ \weier $ and $ \edwards $.
	Moreover, one extends the definition of $ \alpha $ and $ \beta $  by putting $ \alpha((0, 1)) = \Omega $,
	$ \beta(\Omega) = (0, 1) $, $ \alpha((0, -1)) = (0, 0) $ and $ \beta((0, 0)) = (0, -1) $; and
	$ \beta((t_1, 0)) = \beta((t_2, 0)) = \Omega_1 $, $ \beta((-x_1, \pm s_1)) = \Omega_2 $,
	where $ (t_1, 0), (t_2, 0), (-x_1, \pm s_1) \in \weierof{\K} $, with $ t_1, t_2 \ne 0 $.

	The value $ \Omega_1 $ of the two images $ \mathmbox{\beta((t_1, 0))} $, $ \mathmbox{\beta((t_2, 0))} $
	and the value $ \Omega_2 $ of the two images $ \mathmbox{\beta((-x_1, \pm s_1))} $
	can be directly found by passing to homogeneous
	coordinates. As for the value of $ \beta(\Omega) $, we have that,
	for $ \K = \mathbb{C} $, if $ \wp(z) $ and $ \wp^\prime(z) $ are
	the \keyphrase{Weierstrass elliptic functions} and
	$ P = \roundbrack{\wp(z) + \dfrac{a^\prime}{3}, \dfrac{1}{2} {\wp^\prime(z)}} \in \weierof{\K} $,
	then,
	$ \lim\limits_{z \rightarrow 0} \beta(P) = \lim\limits_{z \rightarrow 0} \roundbrack{\dfrac{2 y_1}{3 x_1} \dfrac{(3 \wp(z) + a^\prime)}{{\wp^\prime(z)}}, \dfrac{3 \wp(z) + a^\prime - 3 x_1}{3 \wp(z) + a^\prime + 3x_1}} = (0,1) = O $, as $ \dfrac{3 \wp(z) + a^\prime}{3 {\wp^\prime(z)}} = o(z) $, that is, $ \beta $ is continuous in $ P = \Omega $.
\end{remark}

\begin{remark}\label{rmk:coherently}
	Since there are two points mapped by $ \beta $ onto $ \Omega_1 $ and
	two points onto $ \Omega_2 $, one sees that it is not possible to coherently
	define $ \alpha(\Omega_1) $ and $ \alpha(\Omega_2) $.
	For this reason the maps $ \alpha $ and $ \beta $ define a
	\keyphrase{birational equivalence} between the two forms. Note
	that $ \weier $ is, indeed, a smooth \keyphrase{projective resolution} of the
	non-smooth curve $ \edwards $.
\end{remark}

\begin{remark}
	Since the map $ \beta $ in {\normalfont (\ref{eq:beta})} transforms
	a line through $ P \in \weierof{\K} $ and $ Q \in \weierof{\K} $ onto the hyperbola through $ \beta(P) $, $ \beta(Q) $,
	$ O^\prime, 2 \Omega_1 $ and $ 2\Omega_2 $, and maps vertical lines onto horizontal lines,
	then $ \beta $ induces a group homomorphism of the corresponding divisor classes group
	({\normalfont cf. \cite{BernsteinLange}} for further details about the group law of Edwards curves).
\end{remark}

\section{The Riemann-Roch space $ \gfspace{D} $ over Edwards curves}\label{sec:2}

In this section, given a divisor $ D \in \mathrm{Div}(\edwards) $, we provide a basis of the Riemann-Roch vector space
\begin{equation*}
	\gfspace{D} = \{f \in {\overline{\K}(\edwards)}^\star : \divp{f} + D \mbox{ is effective}\} \cup \{0\}
\end{equation*}
for an Edwards curve $\edwards$, under the assumption that the support of $ D $ does not contain
the two singular points $ \Omega_1 $ and $ \Omega_2 $.

We recall that a \textit{divisor} is, in this context, an element of the free abelian group $ \mathrm{Div}(\edwards) $
on the points of $ \edwards $, that is, a formal sum $ D = \sum_{n_P \in \mathbb{Z}} n_P P $, with $ P \in \edwardsof{\K} $,
where only finitely many integers $ n_P $ are not zero, and that a \textit{principal divisor} $ D = \divp{g} $ of a function
$ g $ is the sum of the zeros of $ g $ on $ \edwards $ minus the poles of $ g $ on $ \edwards $.
The integer $ \delta = \sum n_P $ is the \textit{degree} of the divisor $ D $ and principal divisors give a subgroup of the
subgroup $ {\mathrm{Div}}^0(\edwards) $ of divisors having degree equal to zero, because any function $ g $ on $ \edwards $
has by Bezout theorem the same number of zeros and poles on $ \edwards $.
The group taken into account is formally the quotient
group $ \dfrac{\operatorname{Div}^0(\edwards)}{\operatorname{Princ}(\edwards)} $.

Also, we recall that any divisor $ D^\prime $ on $ \edwards $ of degree $ k + 1 $, such that
$ \Omega_1 $ and $ \Omega_2 $ do not belong to the support of $ D^\prime $, is linearly equivalent to
$ P + kO $, for a suitable point $ P \in \edwardsof{\K} $ (or $ (k + 1)O $, in the case where $ P = O $), that is,
$ D^\prime = P + k O + \divp{g}$, for a suitable function $ g $. Since the map
\begin{small}
	\begin{equation*}
		\begin{aligned}
			\chi : \gfspace{D^\prime} &\longrightarrow \gfspace{ P + kO }\\
			\gen{F}    &\longmapsto    g \gen{F}
		\end{aligned}
	\end{equation*}
\end{small}
is an isomorphism between $ \gfspace{D^\prime} $ and $ \gfspace{P + kO} $, we confine ourselves to the latter space.

\begin{theorem}\label{thm:rredwards}
	Let $ \edwards $ be an Edwards curve defined, over a field $ \K $ of characteristic different
	from $ 2 $, by the equation $ x^2 + y^2 = 1 + d x^2 y^2 $,
	where $ d $ is not a square.
	If $ P + kO = D \in \mathrm{Div}(\edwards) $ is a divisor of positive degree $ k + 1 $,
	where $ P = (a, b) $, then $ \mathrm{dim}(\gfspace{D}) = k + 1 $ and
	\begin{equation*}
		\begin{aligned}
			\gfspace{D} =
			\begin{dcases}
				\left< \gen{F_0}, \gen{F_1}, \ldots, \gen{F_k} \right> &\mbox{if } P \ne O\\
				\left< \gen{F_0}, \gen{F_2}, \ldots, \gen{F_{k + 1}} \right> &\mbox{if } P = O
			\end{dcases}
		\end{aligned}
	\end{equation*}
	where $ \gen{F_0}, \gen{F_1}, \ldots, \gen{F_{k + 1}} $ are rational homogeneous functions defined as follows:
	\begin{equation*}
		\begin{aligned}
			\gen{F_0} &= \dfrac{Z}{Z}\\
			\gen{F_1} &=
			\begin{dcases}
				\begin{aligned}
					&\dfrac{Z}{X} &\mbox{if } P &= O^\prime = (0, -1)\\
				 	&\dfrac{(X + Z)(Y + Z)}{X Y} &\mbox{if } P &= H = (1, 0)\\
					&\dfrac{(X - Z)(Y + Z)}{X Y} &\mbox{if } P &= H^\prime = (-1, 0)\\
					&\dfrac{(Y + b Z) \cdot X}{(X - a Z) \cdot (Y - Z)} &\mbox{if } P &\notin \{O^\prime, H, H^\prime\}\\
				\end{aligned}
			\end{dcases}\\
			\gen{F_i} &=
			\begin{dcases}
				\begin{aligned}
					&\dfrac{Z^h}{{(Y - Z)}^h} &\mbox{if } i &= 2h\\
					&\dfrac{(Y + Z) Z^h}{X {(Y - Z)}^h} &\mbox{if } i &= 2h + 1\\
				\end{aligned}
			\end{dcases}
		\end{aligned}
	\end{equation*}
	for $ 2 \le i \le k + 1 $.
\end{theorem}

\begin{proof}
	Since $ P $ is different from $ \Omega_1 $ and $ \Omega_2 $, we can take the
	point $ \alpha(P) \in \weierof{\K} $, where $ \mathmbox{\alpha: \edwardsof{\K} \longrightarrow\weierof{\K}} $ is the map
	defined in (\ref{eq:alpha}). Hence, the (surjective) map
	\begin{equation*}
		\mathmbox{\alpha^{-1} = \beta: \weierof{\K}\longrightarrow\edwardsof{\K}}
	\end{equation*}
	in (\ref{eq:beta}), induces an (injective) homomorphism
	$ g \mapsto g\circ\beta $ from $ \gfspace{P+kO} $ to
	$ \gfspace{\alpha(P) + k\Omega} $,
	because $ \mathmbox{\beta\roundbrack{\divp{g\circ\beta}} = \divp{g}} $
	for any function $ g \in \gfspace{P + kO}$.

    Since $ \weier $ is smooth, by the formula of Riemann-Roch,
    the dimension of $ \gfspace{\alpha(P) + k\Omega} $ is $ k + 1 $, and we are left with exhibiting
    $ k + 1 $ linearly independent functions in $ \gfspace{P + kO} $, as manifestly $ \gfspace{P + (i - 1)O} $
    is contained in $ \gfspace{P + iO} $, for $ i = 1, \ldots, k + 1 $.

	\noindent For $ i = 0 $ the assertion follows, because $ \divp{\dfrac{Z}{Z}} = 0 $ and for every
	$ P \in \edwardsof{\K} $ we have that $ \divp{\gen{F_0}} + P $ is effective.

	\begin{figure}[H]
		\centering
		\includegraphics[width=0.4\textwidth, keepaspectratio]{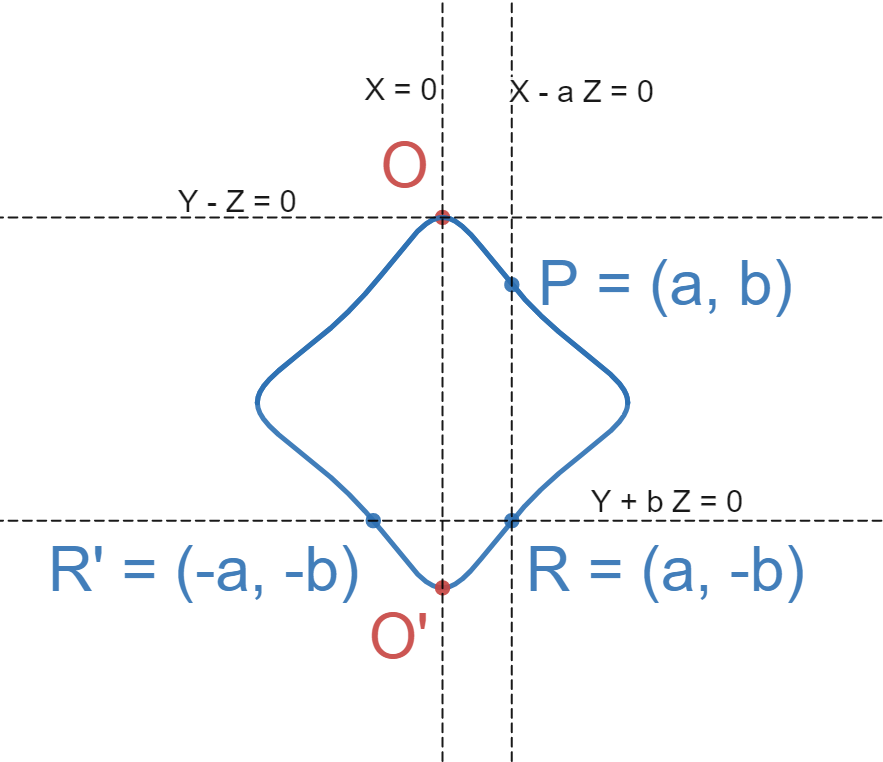}%
		\caption{Edwards real curve with $ d = -8 $. The points $ P, R, R^\prime $ have
			the following coordinates $ (a, b) $, $ (a, -b) $, $ (-a, -b) $, respectively}%
		\label{fig:edwards_proof}%
	\end{figure}

	\noindent Recalling that $ O^\prime = (0, -1) $, $ H = (1, 0) $, $ H^\prime = (-1, 0) $,
	and putting $ R = (a, -b) $ and $ R^\prime = (-a, -b) $ for $ P = (a, b) $ as in \figref{edwards_proof},
	for $ i = 1 $ we have that:
	\begin{equation*}
		\begin{aligned}
			&\begin{aligned}
				\divp{\dfrac{Z}{X}} &= 2H^\prime + 2\Omega_2 + 2O^\prime + 2\Omega_1) - (H + H^\prime + 2\Omega_1 + O + O^\prime + 2\Omega_2) \\ &= H^\prime + O^\prime - H - O,
			\end{aligned}\\
			&\begin{aligned}
				\divp{\dfrac{(X + Z)(Y + Z)}{X Y}} &= (2H^\prime + 2\Omega_2 + 2O^\prime + 2\Omega_1) -
				(H + H^\prime + 2\Omega_1 + O + O^\prime + 2\Omega_2) = \\ &= H^\prime + O^\prime - H - O,
			\end{aligned}\\
			&\begin{aligned}
				\divp{\dfrac{(X - Z)(Y + Z)}{X Y}} &= (2H + 2\Omega_2 + 2O^\prime + 2\Omega_1) -
				(H + H^\prime + 2\Omega_1 + O + O^\prime + 2\Omega_2) =\\&= H + O^\prime - H^\prime - O,
			\end{aligned}\\
			&\begin{aligned}
				\divp{\dfrac{(Y + b Z ) \cdot X}{(X - a Z) \cdot (Y - Z)}} &=
				(R + R^\prime + 2 \Omega_1 + O + O^\prime + 2 \Omega_2) - (P + R + 2 \Omega_2 + 2 O + 2 \Omega_1) = \\
				& = R^\prime + O^\prime - P - O.
			\end{aligned}\\
		\end{aligned}
	\end{equation*}
	Hence, we have that $ \divp{\gen{F_1}} + P + O $ is effective for any suitable $ P $.

	\noindent Additionally, for $ i \ge 2 $, we have that:
	\begin{equation*}
		\begin{aligned}
			&\begin{aligned}
				\divp{\dfrac{Z^h}{{(Y - Z)}^h}} = (2h \Omega_1 + 2h \Omega_2) - (2h O + 2h \Omega_1) = 2h \Omega_2 - 2h O,
			\end{aligned}\\
			&\begin{aligned}
				\divp{\dfrac{(Y + Z) Z^h}{X {(Y - Z)}^h}} &= (2 O^\prime + 2 \Omega_1 + 2h \Omega_1 + 2h \Omega_2) -
				(O + O^\prime + 2 \Omega_2 + 2h O + 2h \Omega_1) = \\
				&= O^\prime + 2 \Omega_1 + (2h - 2) \Omega_2 - (2h + 1) O,
			\end{aligned}
		\end{aligned}
	\end{equation*}
	hence, $ \divp{\gen{F_i}} + P + i O $ is effective in both the cases $ i = 2h $ and $ i = 2h + 1 $.

	So, every function $ \gen{F_i} $ is such that $ \divp{\gen{F_i}} + D $ is
	effective if $ D = P + kO $ or $ D = (k + 1)O $.
	In order to complete the proof, it is necessary to show that all these functions are linearly independent,
	but this follows from standard, elementary, arguments of linear algebra.

	We note that in the case $ D = (k + 1)O $ we simply remove $ \gen{F_1} $ and we add
	$ \gen{F_{k+1}} $, thus, also in this case, we have $ k + 1 $ linearly
	independent functions.
\end{proof}

\begin{remark}
	We note that it is not possible to extend the proof about $ \gfspace{D} $ in
	{\normalfont \thmref{rredwards}} when $ P $ is equal to $ \Omega_1 $ or $ \Omega_2 $, because
	the map $ \beta $ is not invertible on these points {\normalfont (see \remarkref{coherently})}.
	Riemann-Roch spaces on curves having singular points are the subject of {\normalfont $ \S $ IV.2} in
	{\normalfont \cite{Serre1988}}.
\end{remark}

\subsection{Computational cost}\label{subsec:1}

Recalling that the costs of modular addition, multiplication, and inversion over $ \GF{q} $
are $ \ord{\ln(q)} $, $ \ord{\ln^2(q)} $, $ \ord{\ln^3(q)} $, respectively,
we now compute the cost of evaluating at a point $ P \in \edwardsof{\GF{q}} $ each element of the basis of
$ \gfspace{D} $.

We firstly note that we can compute $ \gen{F_{2h}} $ from $ \gen{F_{2h - 2}} $, and $ \gen{F_{2h + 1}} $ from
$ \gen{F_{2h}} $, as
\begin{equation*}
	\begin{dcases}
		\begin{aligned}
			\gen{F_{2h}} &= \gen{F_2} \gen{F_{2h - 2}} &\mbox{if } h &\ge 2,\\
			\gen{F_{2h + 1}} &= \dfrac{Y + Z}{X} \gen{F_{2h}} &\mbox{if } h &\ge 1,
		\end{aligned}
	\end{dcases}
\end{equation*}
that is, at each step we have to perform a single multiplication times the last
(or the second-last) value.
Moreover, we can pre-calculate the value of the function $ \dfrac{Y + Z}{X} $ at $ P $
with a cost $ \cost{\dfrac{Y + Z}{X}} = \ord{\ln(q)} + \ord{\ln^2(q)} + \ord{\ln^3(q)} \approx \ord{\ln^3(q)} $
to further speed up the computation.

Therefore, the maximal global cost $ \cost{\gfspace{D}} $ of
evaluating the first $ k $ functions $ \gen{F_i} $
is $ \cost{\gfspace{D}} = \ord{k \cdot \ln^2(q) + \ln^3(q)} $.

\section{AG Goppa codes on Edwards curves}\label{sec:3}

In this section, we construct the generating matrix and the parity-check matrix for a $ {[n, k, d]}_q $ AG Goppa code
for an Edwards curve $ \edwards $ over $ \GF{q} $, compute the computational cost, and give a small example.

\subsection{Goppa code for an Edwards curve} \label{subsec:notations}

In the following, we adapt the definition of a Goppa code to our case.

\begin{definition}\label{def:aggoppa}
	Let $ D = P + (k - 1)O $ be a divisor of positive degree $ \delta D = k $ of
	the Edwards curve $ \edwards $ over $ \GF{q} $, where $ q = p^t $ and $ p $ is an odd prime number.
	Let $ \gfspace{D} $ be the Riemann-Roch space, let $ T = \{P_1, \ldots, P_n\} $ be
	a set of $ n > k - 1 $ points such that, for $ i = 1, \ldots, k $ and $ j = 1, \ldots, n $,
	$ G_{ij} = \gen{F_{i-1}}(P_j) $, where $ \{\gen{F_i}\} $ is a basis of $ \gfspace{D} $,
	$ P_j \in \edwards $, and $ P_j \notin \mathrm{supp}(D) $.
	Let $ (G_{ij}) = G \in {\GF{q}}^{k \times n} $ be the $ k \times n $ matrix,
	we define the $ {[n, k, d]}_q $ \keyphrase{AG Goppa code}
	$ \mathcal{C}_G = \{\nvector{c} \in {\GF{q}}^n : \nvector{c} = \nvector{a} \cdot G, \nvector{a} \in {\GF{q}}^k \} $.
\end{definition}

\begin{remark}
	We note that $ G $ is well defined because all points $ P_j \in T $ do not belong to the support
	of $ D $ which contains the poles of each $ \gen{F_i} $.
\end{remark}

\begin{theorem}
	If $ \mathcal{C}_G $ is the AG Goppa code of \defref{aggoppa}, then the
	minimum distance $ d $ of this code is such that $ d \ge n - \delta D = n - k $.
\end{theorem}

\begin{proof}
	It follows from the same, classic, proof of AG Goppa codes over curves.
\end{proof}

If we order the points in $ T $ so that the first $ k $ columns of the generating matrix
\begin{equation*}
	G = \begin{pmatrix}
		\gen{F_0}(P_1) & \gen{F_0}(P_2) & \cdots & \gen{F_0}(P_n)\\
		\gen{F_1}(P_1) & \gen{F_1}(P_2) & \cdots & \gen{F_1}(P_n)\\
		\vdots & \vdots & \ddots & \vdots \\
		\gen{F_{k - 1}}(P_1) & \gen{F_{k - 1}}(P_2) & \cdots & \gen{F_{k - 1}}(P_n)\\
	\end{pmatrix} \in {\GF{q}}^{k \times n}
\end{equation*}
of the Goppa code $ \mathcal{C}_G $ are linearly independent, e.g. by applying the Gauss-Jordan method,
then $ G $ can be reduced in its \keyphrase{standard form} $ [I_k \vert M] $,
where $ I_k $ is the identity matrix of order $ k $ and $ M \in {\GF{q}}^{k \times (n - k)} $.
Once $ G $ is in standard form, the parity-check matrix $ H \in {\GF{q}}^{(n - k) \times n} $
of this Goppa code, that is, the matrix such that $ G \cdot H^T = \nvector{0} $ and $ H \cdot {\nvector{y}}^T = \nvector{0} $
for every code word $ \nvector{y} \in \mathcal{C}_G $, is simply $ H = [-M^T \vert I_{n - k}] $.
Thus, the code $ \mathcal{C}_G $ is also defined as $ \{\nvector{y} \in {\GF{q}}^n : H \cdot {\nvector{y}}^T = \nvector{0}\} $.

\subsection{Computational cost of constructing a Goppa code}\label{subsec:2}

In order to compute the generating matrix $ G $ we need to evaluate each of the $ n $ points
in the set $ T $ for each element of the basis of $ \gfspace{D} $, that is, we have a computational cost of
$ \ord{n \cdot \cost{\gfspace{D}}} $ because $ G $ is a matrix of size $ k \times n $.
Moreover, the cost of computing the parity-check matrix depends on the method used to solve
the linear system $ G \cdot \nvector{x} = 0 $.  For instance,
if we used the \keyphrase{Gauss-Jordan method} to reduce the matrix
$ G $ to its standard form, then the cost would be $ \ord{\max{(n, k)}^3} = \ord{n^3} $.

Hence, the global computational cost of constructing a Goppa code over $ \edwardsof{\GF{q}} $ is:
\begin{equation*}
	\ord{n \cdot \cost{\gfspace{D}}} + \ord{\max{(n, k)}^3} = \ord{n \cdot (k \cdot \ln^2(q) + \ln^3(q))} + \ord{n^3}.
\end{equation*}

In particular, for $ k \cdot \ln^2(q) + \ln^3(q) < n^2 $, the computational cost is $ \ord{n^3} $.
However, if $ k \cdot \ln^2(q) + \ln^3(q) > n^2 $,
for instance if we were working with very large finite fields ($ q \gg 2 $),
the overall computational cost would be $ \ord{n \cdot (k \cdot \ln^2(q) + \ln^3(q))} $.

\subsection{A small example}\label{subsec:goppa_edwards}

Let $ \K = \mathrm{GF}(17) $ and let $ \edwards $ be the Edwards
curve defined by the equation $ x^2 + y^2 = 1 + 10 x^2 y^2 $.
There are $ 24 $ affine points on this curve:
\begin{equation*}
	\begin{aligned}
		\{&(0, 1), (0, 16), (1, 0), (2, 2), (2, 15), (3, 6), (3, 11), (5, 8), (5, 9), (6, 3), (6, 14), (8, 5), \\
		  &(8, 12), (9, 5), (9, 12), (11, 3), (11, 14), (12, 8), (12, 9), (14, 6), (14, 11), (15, 2), (15, 15), (16, 0)\},
	\end{aligned}
\end{equation*}
and the two points at infinity.

Let $ D = (2, 15) + 4O $ be the divisor defining $ \gfspace{D} $, so the degree $ k $ is $ 5 $. Let
$ T = \{P_1 = (5, 8), P_2 = (5, 9), P_3 = (6, 3), P_4 = (6, 14), P_5 = (8, 5), P_6 = (8, 12), P_7 = (9, 5)\} $
be the set of points such that $ P_j \notin \operatorname{\mathrm{supp}}(D) $ defining
the generating matrix $ G $ of $ \mathcal{C}_G $, thus $ n = 7 $.
Applying \thmref{rredwards}, the vector space $ \gfspace{D} $ has the following basis:
\begin{equation*}
	\begin{aligned}
		\gfspace{D} &= \left<\gen{F_0}, \gen{F_1}, \gen{F_2}, \gen{F_3}, \gen{F_4}\right> =\\
		&= \left<1, \dfrac{x (y + 15)}{(x - 2)(y - 1)}, \dfrac{1}{y - 1}, \dfrac{y + 1}{x (y - 1)}, \dfrac{1}{{(y - 1)}^2}\right>,
	\end{aligned}
\end{equation*}
whereas the generating matrix $ G = (G_{ij})$ of $ \mathcal{C}_G $ is
defined by putting $ G_{ij} = \gen{F_{i-1}}(P_j) $, that is,
\begin{equation*}
	\begin{aligned}
		G = \begin{pmatrix}
			1 & 1 & 1 & 1 & 1 & 1 & 1\\
			16 & 5 & 5 & 4 & 1 & 11 & 4\\
			5 & 15 & 9 & 4 & 13 & 14 & 13\\
			9 & 13 & 6 & 10 & 14 & 10 & 3\\
			8 & 4 & 13 & 16 & 16 & 9 & 16\\
		\end{pmatrix}.
	\end{aligned}
\end{equation*}
We now compute the parity-check matrix $ H $
by solving the linear system $ G \cdot \nvector{x} = \nvector{0} $, which reduces to:
\begin{equation*}
	\begin{dcases}
		\begin{aligned}
			x_1 &= 7 x_6 + 2 x_7\\
			x_2 &= 3 x_6 + 12 x_7\\
			x_3 &= x_6 + 9 x_7\\
			x_4 &= 13 x_6 + 12 x_7\\
			x_5 &= 9 x_6 + 15 x_7\\
		\end{aligned}
	\end{dcases},
\end{equation*}
that is,
\begin{gather}
	H = \left(\begin{array}{ccccc|cc}
		7 & 3 & 1 & 13 & 9 & 1 & 0\\
		2 & 12 & 9 & 12 & 15 & 0 & 1\\
	\end{array}\right).
\end{gather}
Finally, the minimum distance for this code is $ d \ge n - \delta D = 7 - 5 = 2 $. Moreover,
we know from the Singleton theorem that $ d \le n - k + 1 $ for a $ {[n, k, d]}_q $ code,
that is, $ d \le 7 - 5 + 1 = 3 $ as $ k = \delta D = 5 $.
Hence, $ d = 2 $ or $ d = 3 $, but one can easily check that $ 2 $ columns of $ H $ are always
linearly independent, so $ d = 3 $.
Therefore, we have a $ {[7, 5, 3]}_{17} $-Goppa \keyphrase{MDS} code.

\section{Conclusions}

Edwards curves have been recently introduced for their applications in cryptography.
In this paper, we provided a basis for the Riemann-Roch space of a divisor on these curves, and we used
this basis for the construction of the generating matrices of the AG Goppa codes, thus
providing a possible application of Edwards curves to Coding theory, as well.

\normalem
\printbibliography[title={References}]

\end{document}